\documentclass[a4paper,reqno]{amsart}
\usepackage{amssymb}
\usepackage{latexsym}
\usepackage{amsmath}
\usepackage{euscript}
\usepackage{graphics,color}
\usepackage[all]{xy}
\usepackage[margin=3cm]{geometry}
\usepackage{MnSymbol} 

\newcommand{\be}{\begin{equation}}
\newcommand{\ee}{\end{equation}}
\newcommand{\ba}{\begin{eqnarray}}
\newcommand{\ea}{\end{eqnarray}}
\newcommand{\baa}{\begin{eqnarray*}}
\newcommand{\eaa}{\end{eqnarray*}}
\newcommand{\bb}{\begin{thebibliography}}
\newcommand{\eb}{\end{thebibliography}}

\newcommand{\bi}[1]{\bibitem{#1}}
\newcommand{\lab}[1]{\label{#1}}
\newcommand{\re}[1]{(\ref{#1})}

\newcounter{my}
\newcommand{\he}%
   {\stepcounter{equation}\setcounter{my}%
   {\value{equation}}\setcounter{equation}0%
   }%
\newcommand{\she}%
   {\setcounter{equation}{\value{my}}%
    }%

\renewcommand\t{\tilde}

\newcommand\vphi{\varphi}

\newcommand\ve{\varepsilon}

\newtheorem{pr}{Proposition}

\theoremstyle{definition}

\numberwithin{equation}{section}

\begin{document}

\title{The Heun operator of Hahn type}

\author{Luc Vinet}
\author{Alexei Zhedanov}

\address{Centre de recherches \\ math\'ematiques,
Universit\'e de Montr\'eal, P.O. Box 6128, Centre-ville Station,
Montr\'eal (Qu\'ebec), H3C 3J7}

\address{Department of Mathematics, School of Information, Renmin University of China, Beijing 100872,CHINA}

\begin{abstract}
The Heun-Hahn operator on the uniform grid is defined. 	This operator is shown to map polynomials of degree $n$ to polynomials of degree $n+1$, to be tridiagonal in bases made out of either Pochhammer or Hahn polynomials and to be bilinear in the operators of the Hahn algebra. The extension of this algebra that includes the Heun-Hahn operator as generator is described. Biorthogonal rational functions on uniform grids are shown to be related to this framework.
\end{abstract}

\keywords{}


\maketitle

\section{Introduction}
\setcounter{equation}{0}
The purpose of this  paper is to introduce and characterize the difference 
analog of the Heun operator which we will call the Heun operator of the 
Hahn type. The starting point  will be the construction of the most general 
second order difference operator on the uniform grid that maps polynomials of 
degree $n$ to polynomials of degree $n+1$.  This parallels the similar 
property  that the Heun operator possesses in the continuum. It will be 
observed that the operator thus constructed is tridiagonal in bases either 
made of Pochhammer polynomials or of Hahn polynomials. Finally, the 
connection with the concept of algebraic Heun operator will be established. 

The notion of algebraic Heun operator associated to bispectral problems was introduced in
recent work \cite{GVZ_band}. Essentially, this operator consists in the generic bilinear combination of the
two operators determining the bispectral framework. The name is motivated by
the fact that the operator thus constructed by considering the Jacobi polynomials and using the
hypergeometric operator and multiplication by the variable, precisely coincides with the standard
Heun operator \cite{GVZ_Heun}.

All polynomials of the Askey scheme  \cite{KLS} are known to offer bispectral
situations. It is therefore possible following \cite{GVZ_band} to associate an algebraic
Heun operator to each family of the scheme. In this context, the definition
of the algebraic Heun operator is intimately connected to a procedure referred to as tridiagonalization
\cite{IK1}, \cite{IK2}, \cite{GIVZ} which has been used to develop a bottom-up theoretical description
of orthogonal polynomials where families with higher number of parameters are constructed and
characterized from those with a lesser number of parameters.

It is by now well appreciated that the bispectral properties of the polynomials of the Askey
scheme can be encoded in algebras of quadratic type that are referred to as the Racah and Askey-Wilson
algebras in the case of the polynomials at the top of the hierarchies \cite{Zh_AW}, \cite{GLZ_Annals}. 
The algebras associated to the other families are similarly identified by the names of
the corresponding polynomials. Algebraically, the effect of tridiagonalization is hence to embed the more complicated algebras
into the simpler ones. One can thus insert the Racah algebra into the Jacobi one. It should be
stressed here 
that the most general tridiagonalization that provides the corresponding algebraic Heun
operator transcends the framework of orthogonal polynomials.  In this connection we shall offer a description of the algebraic structure that accompanies this extension in the Heun-Hahn case. This will provide the simplest possible generalization of the Heun - Jacobi picture which arise from the complete tridiagonalization of the hypergeometric operator \cite{GVZ_Heun}.  It will also be observed that certain biorthogonal rational functions arise as solutions of the generalized eigenvalue problems associated to Heun-Hahn difference operators.

The outline of the paper is as follows. In Section 2, we will briefly review the properties of the standard differential Heun operator whose correspondence in the discrete domain we will aim to obtain. This task will be performed in Section 3. We shall recall in Section 4 how the bispectrality of the Hahn polynomials is accounted algebraically and shall introduce the Hahn algebra. How the Heun-Hahn operator is constructed as the algebraic Heun operator associated to this algebra will be shown in Section 5. Section 6 will discuss the generalization of this algebra when the Heun-Hahn operator is taken as a generator. Section 7 will be devoted to certain generalized eigenvalue problems associated to the difference Heun-Hahn operator and to biorthogonal rational functions that arise as solutions. Concluding remarks will be found in Section 8.

\section{Basic properties of the standard Heun operator}
\setcounter{equation}{0}
The ordinary Heun operator has the expression (up to an affine transformation $x \to \alpha x +\beta$ of the independent variable) \cite{Ronveaux}
\be
W = \pi_3(x) \partial_x^2 + \pi_2(x) \partial_x + \pi_1(x)\mathcal{I}, \lab{ord_Heun} \ee
where $\pi_i(x)$ are arbitrary polynomials of degrees $\le i$ and where $\mathcal{I}$ is the identity operator. We will assume additionally that $\deg(\pi_3(x)) \ge 1$. This means that we will exclude the case when $\pi_3(x)$ is a constant.

Simple considerations show that this operator can be characterized by one of the two equivalent conditions:

{\bf (i)} It is the most general second order differential operator which sends any polynomial of degree $n$ to a polynomial of degree $n+1$.

{\bf (ii)} It is the most general second order differential operator which is tridiagonal with respect to the "monomial" basis $\varphi_n(x)=(x-\ve)^n, \: n=0,1,2,\dots$ with some parameter $\ve$ depending on roots of the polynomial $\pi_3(x)$. In particular, if one of these roots is zero (i.e. $\pi_3(0)=0$) then  $\varphi_n(x)=x^n$. The term "tridiagonal" means that
\be
W \varphi_n(x) = \{\vphi_{n+1}(x), \vphi_n(x), \vphi_{n-1}(x)\}, \lab{3-diag_ord_phi} \ee
(i.e. $W \vphi_n(x)$ is a linear combinations of the terms in the symbol $\{\dots\}$).

\vspace{3mm}

Indeed, applying successively the second order differential operator 
\be
W = A(x) \partial_x^2 + B(x) \partial_x + C(x) \mathcal{I} \lab{arb_W_der} \ee
to the constant, linear and quadratic monomials $1,x,x^2$ and demanding that the resulting polynomials  have degrees 1,2 and 3, we arrive at the operator having expression \re{ord_Heun}. Obviously this operator is the most general operator with property (i). Property (ii) follows from the simple observation that by a shift of the argument $x \to x+c$, one can always achieve $A(0)=0$ and hence $A(x)=x \pi_2(x)$ with $\pi_2(x)$ a polynomial of degree $\le 2$. Then property (ii) is obvious with respect to the monomial basis $x^n$. The inverse statement can easily be verified by taking $n=0,1,2$ in \re{3-diag_ord_phi}.

\vspace{3mm}

In \cite{GVZ_Heun} it was showed that there are two more equivalent conditions for the operator $W$ to be the Heun operator:

{\bf (iii)} The operator $W$ is an arbitrary bilinear combination
\be
W= \tau_1 XY + \tau_2 YX + \tau_3 X + \tau_4 Y + \tau_0 \mathcal{I} \lab{W-bil_ord} \ee
of the two operators $X$ and $Y$, where $X$ is the multiplication by $x$
\be
X f(x) = xf(x) \lab{X_mul} \ee
and $Y$ is the hypergeometric operator
\be
Y = \pi_2(x) \partial_x + \pi_1(x) \lab{Gauss_gen} \ee
with polynomials $\pi_2(x), \pi_1(x)$ of degrees 2 and 1.

{\bf (iv)} Up to an affine transformation of the argument $x$, the operator $W$ is the most general second order differential operator which is tridiagonal with respect to the Jacobi polynomials $Q_n^{(\alpha, \beta)}(x)$:
\be
W Q_n(x) = \{Q_{n+1}(x), Q_n(x), Q_{n-1}(x)\}. \lab{Jac_W_3diag}
\ee

\vspace{3mm}

The main purpose of the present paper is to introduce an analog of the Heun operator $W$ on the finite uniform grid. We will show that similarly to the ordinary Heun operator, the discrete Heun operator can be characterized by 4 equivalent conditions very similar to the above conditions (i)-(iv).

\section{Construction of the Heun operator on the uniform grid}
\setcounter{equation}{0}

Up to an affine transformation, any finite uniform grid on the real line can be reduced to the grid $x=0,1,\dots, N$ with some positive integer $N$. The {\it difference} Heun operator $W$ can thus be characterized by one of the following two equivalent conditions:

\vspace{2mm}

{\bf (i)} The operator $W$ is the most general second order difference operator on the uniform grid $x=0,1,\dots, N$ which sends any polynomial of degree $n$  to a polynomial of degree $n+1$ (it is assumed that $0 \le n \le N$).

\vspace{2mm}

{\bf (ii)} It is the most general second order differential operator which is tridiagonal with respect to the Pochhammer basis $\varphi_n(x)=(-x)_n, \: n=0,1,2,\dots,N$,
where the Pochhammer symbol is defined as usual by $(x)_n=x(x+1) \dots, (x+n-1)$:
\be
W \varphi_n(x) = \{\vphi_{n+1}(x), \vphi_n(x), \vphi_{n-1}(x)\}. \lab{3-diag_ord_Pochh} \ee

\vspace{2mm}

Let us start with the  construction of the difference Heun operator in accordance with the statement (i). Assume that $W$ is a generic second order difference operator acting on the uniform grid:
\be
W = A_1(x) T^+ + A_2(x) T^- + A_0(x) \mathcal{I} 
\lab{W_gen} \ee
with some functions $A_i(x), i=0,1,2$. Here we use the standard notation for the shift operators,
\be
T^{\pm} f(x) = f(x \pm 1). \lab{T_pm} \ee
Act first with the operator $W$ on the constant:
\be
W\{1\} = A_1(x) +A_2(x) + A_0(x). \lab{W_const} \ee 
By condition (i) we have from \re{W_const}
\be
A_0(x) = \pi_1(x) -A_1(x) - A_2(x), \lab{A0_12} \ee
where $\pi_1(x)$ is a polynomial of first degree. Applying the operator $W$ to the function $x$,  we have similarly
\be
A_1(x) (x+1) + A_2(x) (x-1) + A_0(x) x = \pi_2(x) \lab{2c_A12} \ee
with some polynomial $\pi_2(x)$ of second degree. Taking into account \re{A0_12}, we get
\be
A_2(x)-A_1(x) = \tilde \pi_2(x) = x \pi_1(x) - \pi_2(x), \lab{A2-A1} \ee
where $\deg(\t \pi_2(x)) \le 2$. Finally, applying the operator $W$ to the monomial $x^2$ we find
\be
A_1(x) (x+1)^2 + A_2(x) (x-1)^2 + A_0(x) x^2 = \pi_3(x) \lab{3c_A12} \ee
with $ \pi_3(x)$ a polynomial of maximal degree equal to 3. Again, from \re{A0_12} and \re{A2-A1}, we obtain
\be
A_1(x) + A_2(x) = \tilde \pi_3(x) \lab{A1+A2} \ee
with $\t \pi_3(x)$ also a polynomial of third degree. 
It follows from these results that $A_0(x), A_2(x), A_3(x)$ have the following general expressions
\be
A_0(x) = \pi_1(x) -\t \pi_3(x), \; A_1(x) = \frac{\t \pi_3(x) - \t \pi_2(x)}{2}, \; A_2(x) = \frac{\t \pi_3(x) + \t \pi_2(x)}{2}. \lab{A012} \ee
We shall only consider here the generic case where $\deg(\pi_3(x)) =3$. Then, all the coefficients $A_i(x), \: i=0,1,2$ are third degree polynomials and the leading coefficients of $A_1(x)$ and $A_2(x)$ are the same.

Moreover, the operator $W$ should act on the finite grid $0,1,\dots,N$. This is possible only if the coefficient $A_1(x)$  has $x-N$ as a factor while the coefficient $A_2(x)$ has $x$ as a factor.    

We thus arrive at the expressions
\be
A_1(x) = (x-N) (\kappa x^2 + \mu_1 x + \mu_0), \quad A_2(x) = x (\kappa x^2 + \nu_1 x + \nu_0), \quad A_0(x) = -A_1(x) - A_2(x) + r_1 x+ r_0 \lab{expl_gen_A} \ee
with 7 arbitrary parameters $\kappa, \mu_i, \nu_i, r_i, i=0,1$. 

It remains to show that the operator $W$ thus constructed, satisfy the property (i) for all $n=3,4,\dots, N$. It is directly verified that 
\be
W \{x^n\} = \sigma_n^{(1)} x^{n+1} + O(x^n), \quad n=3,4,\dots, N, \lab{leading_sigma} 
\ee
where
\be
\sigma_n^{(1)} = \left(\mu_1-\nu_1 +\kappa(n-1-N)\right)n + r_1. \lab{sigma} \ee
From this formula it follows that the operator $W$ sends any polynomial of degree $n=0,1, \dots, N$ to a polynomial of degree $n+1$ (provided $r_1 \ne 0$). This proves the statement (i).

We can now establish the property (ii). Indeed, it is directly verified that for the polynomial basis $\vphi_n(x)=(-1)^n(-x)_n$ we have
\be
W \vphi_n(x) = \sigma_n^{(1)} \vphi_{n+1}(x) + \sigma_n^{(2)} \vphi_n(x) + \sigma_n^{(3)} \vphi_{n-1}(x) \lab{W_phi_Pochh} \ee
where
\be
\sigma^{(2)}_n =r_0 +(\mu_0 + r_1 - \nu_0 -\mu_1 N-\nu_1 n)n  +n(2n-1)(\mu_1-\kappa(N-n+1))
\ee
and
\be
\sigma^{(3)}_n = -n(N-n+1) \left(\mu_0 +(n-1)(\mu_1 +\kappa(n-1)) \right).   \lab{sigma_3} \ee
The inverse statement is almost obvious.

Let us remark that the difference  Heun operator $W$ possesses several properties which are similar to the properties of the ordinary Heun operator.  In the limit $x \to Nx, \: N \to \infty$,  the operator $W$ becomes the ordinary second-order differential Heun operator. This is obvious from the observation that the Pochhammer basis $\phi_n(x)$ becomes the monomial basis $x^n$ in this limit and that the operator $W$ becomes the differential operator of second order which is 3-diagonal with respect to the monomial basis. This is known to characterize the ordinary Heun operator.

Moreover the operator $W$ admits specific finite-dimensional truncation conditions which lead to "quasi-exactly 
solvable problems" similarly to what has been observed for the the ordinary Heun operator \cite{Tur_Heun}, \cite{Tur_quasi}. Indeed, assume, e.g. that $\xi_{M+1}^{(1)}=0$ for some positive integer $M<N$. Then the difference operator $W$ can be restricted to the linear space of polynomials of degree $\le M$. The eigenvalue problem
\[
W \psi(x) = 0 \lab{W_psi}
\]
can then be exactly solved at least for small $M=1,2,3,4$.

\vspace{3mm}

In Section 5 we shall establish another equivalent characterization of the Hahn-Heun operator, namely:

{\bf (iii)} The operator $W$ is the most general second order difference operator on the uniform grid $x=0,1,\dots, N$ which is tridiagonal with respect to the Hahn polynomials $P_n(x;\alpha, \beta, N)$ i.e.
\be
W P_n(x) = \{ P_{n+1}(x), P_n(x), P_{n-1}(x) \}. \lab{Hahn_W_3diag}
\ee
The monic Hahn polynomials are defined as \cite{KLS}
\be
P_n(x;\alpha,\beta,N) = \kappa_n {_3}F_2\left({-n, n+\alpha+\beta+1, -x \atop \alpha+1, -N} \left | \right . ;1 \right) \lab{Hahn_P} \ee
where
\be
\kappa_n = \frac{(\alpha+1)_n(-N)_n} {(n+\alpha+\beta+1)_n} \lab{kappa_n} 
\ee
is the normalization factor needed to make the polynomials monic
\be
P_n(x)=x^{n} + O(x^{n-1}) .
\ee
Note that by  "most general" we mean that the operator $W$ is defined up to affine transformation $W \to \alpha W + \beta \mathcal{I}$.

\section{The Hahn algebra}
\setcounter{equation}{0}
The quadratic Racah algebra \cite{GZ_preprint}, \cite{GLZ_Annals} consists of 3 generators $K_1,K_2,K_3$ with the commutation relations
\ba
&&[K_1,K_2]=K_3, \nonumber \\
&&[K_2,K_3]= a_1 \{K_1,K_2\} + a_2 K_2^2 + b K_2 +c_1 K_1 + d_1 \mathcal{I}, \nonumber \\
&&[K_3,K_1]= a_1 K_1^2 + a_2 \{K_1,K_2\}+ b K_1 + c_2 K_2 + d_2 \mathcal{I} \lab{RA}
\ea
where $\mathcal{I}$ is the identity operator and $a_1,a_2,\dots, d_1, d_2$ are the structure parameters of the algebra. The Racah algebra naturally describes the eigenvalue problems of the Racah (Wilson) polynomials and arises in various physical and mathematical contexts  where these polynomials appear (see \cite{GZ_preprint}, \cite{GLZ_Annals}, \cite{GVZ_Racah}, \cite{Ter} for more details).

The Hahn algebra is a specialization of the Racah algebra when one of the parameters $a_1$ or $a_2$ becomes zero. Taking for example $a_2=0$, we obtain the relations 
\ba
&&[K_1,K_2]=K_3, \nonumber \\
&&[K_2,K_3]= a \{K_1,K_2\}  + b K_2 +c_1 K_1 + d_1 \mathcal{I}, \nonumber \\
&&[K_3,K_1]= a K_1^2 + b K_1 + c_2 K_2 + d_2  \mathcal{I} . \lab{HA}
\ea
The Hahn algebra describes the eigenvalue problems of the Hahn polynomials (both discrete and continuous)  and provides in particular an interpretation of the Clebsch-Gordan problem \cite{GZ_preprint}, \cite{GLZ_Annals}.

In what follows we concentrate only on the Hahn algebra and its realizations.

There is a "canonical" realization of the Hahn algebra in terms of difference operators. Indeed, let us introduce the operator $X$ which is the multiplication by $x$:
\[
X=K_1 =x \lab{X_x} 
\]
and the operator $Y$ which is the Hahn difference operator \cite{KLS}
\be
Y=K_2 = B(x) T^+ + D(x) T^- -(B(x)+D(x)) \mathcal{I} \lab{Y_def}
\ee
where $T^{\pm}f(x) = f(x\pm 1)$ are shift operators and where
\be
B(x) =(x-N)(x+\alpha+1), \quad D(x) = x(x-\beta-N-1). \lab{BD_Hahn} 
\ee

The Hahn polynomials $P_n(x)$ are eigenfunctions of the  operator $Y$ \cite{KLS}
\be
Y P_n(x) = \lambda_n P_n(x), \lab{Y_P}
\ee
with the eigenvalue
\be
\lambda_n =n(n+\alpha+\beta+1). \lab{eig_H}
\ee
It is directly verified that the operators $X$ and $Y$ (together with their commutator $K_3=[X,Y]$) satisfy the Hahn algebra \re{HA} with
\be
a=-2, \: b=2N+\beta-\alpha, \: c_1=-(\alpha+\beta)(\alpha+\beta+2), \: c_2= -1, \: d_1 =N(\alpha+1)(\alpha+\beta), \: d_2=N(\alpha+1).
\ee

The Hahn algebra encompasses a duality property: in the basis of the Hahn polynomials, the operator $Y$ becomes diagonal \re{Y_P} while the operator $X$ becomes tridiagonal
\[
X P_n(x) = P_{n+1}(x) + b_n P_n + u_n P_{n-1}(x), \lab{X_P_rec} 
\]
where the recurrence coefficients $b_n,u_n$ are given in \cite{KLS}. We will not use their explicit expression in this paper.

\section{The algebraic Heun operator of the Hahn type}
\setcounter{equation}{0}
There is another approach to construct operators of the Heun type. Following the method of \cite{GVZ_band} we can define the {\it algebraic Heun operator} of Hahn type  as the generic bilinear combination of the operators $X,Y$ belonging to the Hahn algebra:
\be
W = \tau_1 XY + \tau_2 YX + \tau_3 X + \tau_4 Y + \tau_0 \mathcal{I} . \lab{W_Heun}
\ee
Note that operators of such type were considered by Nomura and Terwilliger \cite{NT} in the finite-dimensional context of the Leonard pairs as the most general operators which are tridiagonal with respect to dual bases (that diagonalize either operator $X$ or $Y$).

The operators $W$ are generically called algebraic Heun operators because in the special case of the Jacobi algebra (where $Y$ is the Gauss differential hypergeometric operator) the operator $W$ coincides with the ordinary Heun operator (see \cite{GVZ_Heun} for details). Moreover, for other special cases corresponding to Laguerre and Hermite polynomials (where the corresponding algebras degenerate to the Lie algebra $sl_2$ or to the oscillator algebra), the operator $W$ corresponds to the confluent Heun operator. In \cite{GVZ_band}, it was demonstrated that the algebraic Heun operators provide solutions of the generalized Slepian-Pollak-Landau problems for band and time limiting of the operators from the Askey scheme.

 Explicitly, for the generators $X=x$ and $Y$  (see \re{Y_def}) of the Hahn algebra, the  operator $W$ is the second-order difference operator on the uniform grid
\be
W= A_1(x) T^+ + A_2 T^- + A_0(x) \mathcal{I} 
\lab{W_Heun_BDC} \ee
where
\be
A_1(x) = (x-N)(x+\alpha+1)\left( (\tau_1+\tau_2)x + \tau_2+\tau_4\right)
\lab{B1} \ee
\be
A_2(x) = x(x-\beta-N-1)\left( (\tau_1+\tau_2)x + \tau_4-\tau_2\right)
\lab{D1} \ee
\be
A_0(x) = -A_1(x)-A_2(x) + \left( (\alpha+\beta+2)\tau_2 +\tau_3 \right)x +\tau_0 -N(\alpha+1) \tau_2 .
\lab{C1} \ee
Comparison of \re{W_Heun_BDC}-\re{C1} with \re{W_gen} and \re{expl_gen_A} shows that these operators coincide (as was anticipated by the fact that we used the same symbol W). This means that the parameters $\kappa, \mu_i,\nu_i,r_i$ in \re{expl_gen_A} can be expressed in terms of the parameters $\tau_i, \: i=0,1,\dots,4$ and $\alpha, \beta$ and vice versa. For example, there is the simple formula
\be
\kappa= \tau_1 + \tau_2 \lab{kappa_xi} \ee
which shows that the degrees of the polynomials $A_i(x), \: i=0,1,2$ become $\le 2$ (i.e. when $\kappa=0$) if and only if the algebraic Heun operator can be presented as
\be
W= \tau_1 [X,Y] + \tau_3 X + \tau_4 Y + \tau_0 \mathcal{I}. \lab{W-red} \ee
We have thus established the following new characterization property:

\vspace{2mm}

{\bf (iv)} The Hahn-Heun operator is the second-order difference operator on the uniform \\ grid $0,1,\dots,N$ given by the bilinear Ansatz \re{W_Heun} with $X$ and $Y$ the bispectral operators of the Hahn polynomials.

\vspace{2mm}

The definition  \re{W_Heun} has advantages for the algebraic analysis of the Hahn-Heun operators because it separates the parameters of the operator $W$ into two sets: the {\it internal}  parameters $\alpha, \beta$ coming from the Hahn polynomials and the {\it external} parameters $\tau_i, \: i=0,1,\dots,4$ coming from the bilinear combination \re{W_Heun}.

It is now obvious that the operator $W$ is tridiagonal with respect to the Hahn polynomial basis 
\be
W P_n(x) = \xi_{n+1} P_{n+1}(x) +\eta_n P_n(x) +\zeta_n u_n P_{n-1}(x).
\lab{W_3_Hahn} \ee 
The coefficients $\xi_n, \eta_n, \zeta_n$ can easily be calculated using formula \re{X_P_rec} and have expressions similar to those found in the case  of the ordinary Heun operator \cite{GVZ_Heun}
\ba
&&\xi_n = \tau_1 \lambda_{n-1} + \tau_2 \lambda_n + \tau_3, \: \zeta_n=\tau_2 \lambda_{n-1} + \tau_1 \lambda_n + \tau_3, \nonumber \\
&& \eta_n = (\tau_1+\tau_2) \lambda_n b_n + \tau_3 b_n + \tau_4 \lambda_n + \tau_0. \lab{rec_cf_W} 
\ea
The inverse property can also  be checked. Namely, if $W$ is the most general operator which is 3-diagonal with respect to the Hahn polynomials \re{Hahn_W_3diag}, taking $n=0,1,2$, we arrive at coefficients $A_1(x),A_2(x), A_0(x)$ having expressions \re{expl_gen_A}. This means that the Hahn-Heun operator has the additional characterization property that we already mentioned in Section 3:

\vspace{2mm}

{\bf (iii)} The operator $W$ is the most general second order difference operator on the uniform grid $x=0,1,\dots, N$ which is tridiagonal with respect to the Hahn polynomials $P_n(x;\alpha, \beta, N)$:
\be
W P_n(x) = \{ P_{n+1}(x), P_n(x), P_{n-1}(x) \}. \lab{Hahn_W_3diag-1} \ee

Note also that property \re{W_phi_Pochh} follows automatically from the bilinear Ansatz \re{W_Heun} because the operators $X$ and $Y$  are two-diagonal in the basis $\vphi_n(x)$:
\ba
&&X \vphi_n(x) = n \vphi_n(x) - \vphi_{n+1}(x), \nonumber \\ 
&&Y \vphi_n(x) =  n(n+1+\alpha+\beta)\vphi_n(x) + n(N-n+1)(\alpha+n)\vphi_{n-1}(x). \lab{W_phi}
\ea

We thus see that the bilinear Ansatz \re{W_Heun} can explain many features of the difference Heun operator. Summing up, we have established that the Hahn-Heun operator $W$ can be characterized by the four equivalent properties (i)-(iv) similarly to the case of the ordinary Heun operator.

\section{The Heun-Racah algebra and its reduction to the Racah algebra}
\setcounter{equation}{0}
The Heun operator $W$ together with the Hahn operator $Y$ generate an algebra defined by the following {\it cubic}  relations  
\be
[Y,[Y,W]] = g_1 Y^2 + g_2 \{Y,W\} + g_3 Y + g_4 W + g_5 \mathcal{I} \lab{RH_1}
\ee
and 
\be
[W,[W,Y]] = e_1 Y^2 + e_2 Y^3 + g_2 W^2 + g_1 \{Y,W\} + g_3 W + g_6 Y + g_7 \mathcal{I} . \lab{RH_2} 
\ee
In the above relations, the terms with the coefficients $g_i, \: i=1,\dots,7$ have a structure similar to those of the Racah algebra \re{RA}. The two extra-terms $e_1 Y^2$ and  $e_2 Y^3$ are the ones that make  the algebra with relations \re{RH_1} and \re{RH_2} more general than the Racah algebra. The same situation occurs with the commutation relations involving the ordinary Heun and hypergeometric operators, where extra-terms appear with respect to the Racah algebra relations \cite{GVZ_Heun}. It is therefore natural to refer the algebra \re{RH_1} - \re{RH_2} as the Heun-Racah algebra of the Hahn type.

The coefficients of the extra terms (these terms are beyond the scope of the Racah algebra) are
\ba
&&e_2=2(\tau_1+\tau_2)^2, \nonumber \\
&& e_1 = 6 \tau_4^2 +3(\tau_1+\tau_2)(\tau_3 +(2N+\beta-\alpha)\tau_4) - \nonumber \\
&&(\tau_1^2+\tau_2^2)(3N(\alpha+1)-2) -2(3N(\alpha+1)-5) \tau_1 \tau_2 . \lab{c_12} \ea

\vspace{5mm}

It is seen that if (and only if) the conditions
\be
\tau_1+ \tau_2=0 \lab{t21=0}
\ee
and 
\be
\tau_2 \pm \tau_4 =0 \lab{t24=0} 
\ee
are fulfilled, will the extra terms disappear and the algebra \re{RH_1} - \re{RH_2} become the Racah algebra. So, there are two possible cases of such degeneration of the Heun operators corresponding to the choices of sign in \re{t24=0}. We can thus can introduce the two operators $W_1$ and $W_2$. The first operator $W_1$ is
\be
W_1 =\frac{[X,Y]}{2} + \gamma X - \frac{Y}{2} + \ve \mathcal{I}.
\ee
The second operator is
\be
W_2 =-\frac{[X,Y]}{2} - \gamma X - \frac{Y}{2} - \ve \mathcal{I}
\ee
with $\gamma, \ve$ arbitrary parameters. Irrespective of the restrictions that lead to these operators, one can verify directly that the operators $Y,W_1,W_2$ constitute a Racah triple: any pair of these operators satisfy the Racah algebra relations. Moreover, these operators satisfy the condition
\be
Y+W_1+W_2=0 \lab{YWW=0} 
\ee 
from which it follows that the operators $Y, W_1, W_2$ form the equitable Racah algebra \cite{GWH}, \cite{GCZ_equi}.

Explicitly these operators become the first order difference operators
\be
W_1=(x+\alpha+1)(N-x) T^+ + \left(x^2 +((\alpha-\beta)/2 -N+\gamma)x +\ve -N(\alpha+1)/2
 \right) \mathcal{I} \lab{W_1} 
\ee
and 
\be
W_2= x(\beta+N+1-x) T^- + \left({x}^{2}+ \left( -\gamma -N +(\alpha-\beta)/2 \right) x-\ve -N(\alpha+1)/2    \right)\mathcal{I}. \lab{W_2}
\ee
We thus have 
\begin{pr}
The operators $Y,W_1,W_2$ realize equitably the Racah algebra.
\end{pr} 
Let us remark that the defining relations of the algebra  generated by the operator $W$ and the operator $X$ (instead of $Y$) can also be determined. They take the following form:
\be
[X,[X,W]]= e_3 X^3 + g_8 X^2 + g_9 \{X,W\} + g_{10} X + g_{11} W + g_{12} \mathcal{I}    \lab{RXW_1} \ee
and 
\be
[W,[W,X]]= e_4 X^2 + e_5 X^3 + e_6 XWX  + g_9 W^2 + g_8 \{X,W\} +g_{13} X + g_{14} W + g_{15}. \mathcal{I} \lab{RXW_w} \ee
The expressions for the "Racah-type" structure parameters $g_i, i=8, \dots, 15$ and for the extra parameters $e_3,e_4,e_5,e_6$ are rather complicated and will not be recorded here. It is found that these relations only reduce to those of the Racah algebra in the trivial situation where W degenerates to Y (up to affine transformations).

We conclude this section with the following observation. There is a realization of the Hahn algebra with differential operators. Indeed one can check that the pair of the operators
\be
X= x(x-1) \partial_x + q_1(x) \mathcal{I} , \quad Y= x(1-x) \partial_x^2 + t_1(x) \partial_x  \lab{diff_XY} \ee
satisfy the generic Hahn algebra relations (up to affine transformations of the generators $X$ and $Y$). Here $q_1(x)$ and $t_1(x)$ are arbitrary linear functions in $x$ and $Y$ is the ordinary hypergeometric operator.

Realization \re{diff_XY} can be obtained as a special case of the tridiagonalization of the hypergeometric operator $Y$ \cite{GIVZ}. Consider in this context  the algebraic Heun operator of the Hahn type \re{W_Heun} with the only restriction that $\tau_2=-\tau_1$:
\be
W= \tau_1 [X,Y] + \tau_3 X + \tau_4 Y + \tau_0 \mathcal{I}. \lab{W_restr} \ee
Substituting expressions \re{diff_XY} into \re{W_restr}, we find that the operator $W$ has the expression \re{ord_Heun} with arbitrary polynomials $\pi_1(x), \pi_2(x)$ of degrees one and two and with 
\be
\pi_3(x) = x(x-1)(2 \tau_1 x- \tau_1-\tau_4). \lab{pi3_Heun} \ee
In other words, the algebraic Heun operator of Hahn type becomes the ordinary Heun operator. This example shows that the ordinary Heun operator and the Heun operator of the Hahn type are closely related. We should stress, however, that the ordinary Heun operator appears only under the restriction   $\tau_2=-\tau_1$. In the generic case, in the framework of the realization \re{diff_XY} of the operators X and Y, the Heun operator of Hahn type \re{W_Heun} is a differential operator of {\it third} order:
\be
W = -(\tau_1+\tau_2)x^2(x-1)^2 \partial_x^3 + x(x-1) \sigma_1(x) \partial_x^2 + \pi_2(x) \partial_x + \pi_1(x)\mathcal{I} \lab{W_gen_3} \ee
with some polynomials $\sigma_1(x), \pi_1(x)$ of the first degree and a polynomial $\pi_2(x)$ of the second degree.

\section{The difference Heun operator and generalized eigenvalue problems}
\setcounter{equation}{0}
In this section we show that there are special explicit solutions of the difference Heun operator $W$ which are related to the theory of biorthogonal rational functions on the one hand and to the generalized eigenvalue problem  (GEVP) on the other hand.

The monic Hahn-type $R-II$-polynomials (see \cite{Zhe_Pade} for details) are defined by
\be
P_n(x) = \frac{(\alpha)_n (-N)_n}{(\beta+1)_n} \: {_3}F_2 \left( {-n, -x, -\beta-n \atop -N, 1-\alpha-n} ;1\right). \lab{HR} 
\ee 
Introduce the functions 
\be
U_n(x) = (-1)^n \frac{P_n(x)}{(\alpha-x)_n} = (-1)^n \frac{(-N)_n}{(\beta+1)_n} \: {_3}F_2 \left( {-n, -x, \beta+N-n \atop -N, \alpha-x} ;1\right)
\lab{U_def} \ee
which are rational functions of type $[n/n]$ with simple poles at $x=\alpha, \alpha+1, \dots, \alpha+n-1$ (see \cite{Zhe_Pade} for background and nomenclature). These functions are "monic" in the following sense
\[
U_n(x) = 1 + O(x^{-1}) \lab{monic_U}.
\]

These rational functions possess a certain property of {\it bispectrality} in that they satisfy both a three-term recurrence relations and a second-order difference equation on the uniform grid. The relations can be easily derived from the corresponding relations of the Hahn polynomials \cite{KLS} or from the two-term recurrence relations presented in \cite{Zhe_Pade}.

In contrast to the theory of orthogonal polynomials, these relations belong to the class of generalized spectral problems \cite{Zhe_Pade}. As difference equation one has
\be
L_1 U_n(x,n) = \lambda_n L_2 U_n(x,n), \lab{dfr_U} \ee 
where $L_1$ and $L_2$ are the difference operators defined as
\ba
&&L_1 =  \left( -\alpha+x+1 \right) (x-\alpha) \left( x-N \right) T^+ +x (x-\alpha) \left( x+\beta-\alpha-N \right) T^- + \nonumber \\
 &&(x-\alpha) \left( -2\,{x}^{2}+ \left( 2\,\alpha-1+2\,N-\beta \right) x-N \left( 
\alpha-1 \right)  \right) \mathcal{I},   \lab{L1_def} \ea
and 
\be
L_2 =(x-\alpha) \mathcal{I} - x T^- . \lab{L2_def} \ee
The eigenvalue $\lambda_n$ is
\be
\lambda_n = n(N-\beta-n). \lab{lambda} \ee
Both operators $L_1$ and $L_2$ belong to the family of Hahn-Heun operators \re{W_gen} with coefficients $A_1(x), A_2(x), A_0(x)$ given by \re{expl_gen_A}. For instance, the operator $L_2$ corresponds to the choice of the parameters
\be
\kappa=\mu_1=\mu_0=\nu_1=r_1=0, \quad \nu_0=-1, \; r_0=-\alpha .
\lab{par_L2} \ee
This means that the eigenvalue problem can be presented as follows
\be
L(\lambda_n) U_n(x) =0, \lab{lin_L_U} \ee
where 
\be
L(\lambda_n) = L_1 - \lambda_n L_2 \lab{lin_pen_L} \ee
is a linear pencil of two Hahn-Heun operators of a special type.

Note that it was recently observed that some "classical" biorthogonal rational functions satisfy the ordinary Heun equation \cite{MNR}. Our finding can be considered as a possible discrete analog of this result. One can suppose that any "classical" biorthogonal rational functions should similarly satisfy some differential or difference equations of Heun type.

\vspace{5mm}

\section{Conclusion}
\setcounter{equation}{0}
We have constructed a discrete Hahn-Heun operator which possesses the properties (i)-(iv) analogous to those of the ordinary Heun operator. Moreover, the Heun-Hahn operator has the ordinary Heun operator as a limit when $N \to \infty$.

Looking at the algebra formed by the Hahn operator $Y$ and its tridiagonalization $W$, it was found that the corresponding algebra differs from the Racah algebra by two additional terms. This new Heun-Racah algebra of the Hahn type reduces to the Racah algebra for two specializations of the Heun operator. Moreover, taken together, the Hahn operator and the two specialized Heun operators were  seen to yield an equitable presentation of the Racah algebra.

We have showed also that a linear pencil of two Hahn-Heun operators gives a solution of the difference equation for  biorthogonal rational functions on the uniform interpolation grid. All in all, the introduction of the Heun operator of Hahn type has shed rather interesting light on various properties pertaining to the Hahn bispectral problem.

In the future, it would certainly be appropriate to explore the tridiagonalization of the $q$-analogs of the Hahn operator, namely to proceed with the examination of the algebraic Heun operators associated to the little and big $q$-Jacobi polynomials. This is bound to provide q-analogs of the Heun equation. Such equations have already been discussed in \cite{Takemura} and comparing the contents of this reference with the tridiagonalization outcomes should be telling. This study which we plan to undertake soon will also complement the results obtained in \cite{BMVZ} where the little and big $q$-Jacobi polynomials were found to form representation bases for specialized Askey-Wilson algebra and where the little q-Jacobi operator and a tridiagonalization of it were found to realize the equitable embedding of the Askey-Wilson algebra in $U_q(sl(2))$.

\bigskip\bigskip
{\Large\bf Acknowledgments}

We thank Pascal Baseilhac for discussions and Alphonse Magnus for drawing our attention to the preprint \cite{MNR} some time ago.  The research
of L.V. is funded in part by a discovery grant from the Natural Sciences and Engineering Research
Council (NSERC) of Canada. The work of A.Z. is supported by the National
Science Foundation of China (Grant No.11771015).

\vspace{15mm}

\bb{99}

\bi{AC} H. Alnajjar and B. Curtin. {\it Leonard pairs from the equitable
basis of $sl_2$}. Elec. J. Lin. Alg., {\bf 20} :490-–505, 2010.

\bi{BMVZ} P.Baseilhac, X. Martin, L.Vinet and A.Zhedanov, {\it Little  and big q-Jacobi polynomials and  the Askey-Wilson algebra}, arXiv:1806.02656.

\bi{GWH} S. Gao, Y. Wang, and B. Hou. {\it The classification
of Leonard triples of Racah type} . Lin. Alg. Appl.,
{\bf 439}: 1834-–1861, 2013.

\bi{GVZ_Racah} V.Genest, L.Vinet and A.Zhedanov, {\it Superintegrability in Two Dimensions and the Racah-–Wilson Algebra}, Lett.Math.Phys. {\bf 104}, 931--952 (2014). arXiv:1307.5539v1.

\bi{GCZ_equi} V.Genest, L.Vinet and A.Zhedanov, {\it The equitable Racah algebra from three $su(1,1)$ algebras},  J. Phys. A: Math. Theor. {\bf 47} (2014), 025203. Arxiv: 1309.3540.

\bi{GIVZ} V. Genest, M.E.H. Ismail, L. Vinet, A. Zhedanov {\it Tridiagonalization of the hypergeometric operator and the Racah-Wilson algebra}, arXiv:1506.07803.

\bi{GZ_preprint} Ya.I.Granovskii and A.Zhedanov, {\it Exactly solvable problems and their quadratic algebras}, Preprint, DonFTI, 1989.

\bi{GLZ_Annals} Ya. A. Granovskii, I.M. Lutzenko, and A. Zhedanov, {\it Mutual integrability, quadratic algebras,
and dynamical symmetry}. Ann. Phys. {\bf 217} (1992),  1--20.

\bi{GVZ_Heun} F.A.Gr\"unbaum, L.Vinet and A.Zhedanov, {Tridiagonalization and the Heun equation}, J.Math.Physics {\bf 58}, 031703 (2017), arXiv:1602.04840.

\bi{GVZ_band} F.A.Gr\"unbaum, L.Vinet and A.Zhedanov, {Algebraic Heun operator and band-time limiting},  arXiv:1711.07862.

\bi{Hahn_Heun} W. Hahn, {\it On linear geometric difference equations with accessory parameters}, Funkcial.
Ekvac. {\bf 14} (1971), 73–-78.

\bi{IK1} M. E. H. Ismail and E. Koelink. {\it Spectral analysis of certain Schr\"odinger operators}. SIGMA,
{\bf 8}: 61–-79, 2012.

\bi{IK2} M. E. H. Ismail and E. Koelink. {\it The J-matrix method}.
Adv.Appl.Math., {\bf 56}, 379--395, 2011.

\bibitem{KLS} R. Koekoek, P.A. Lesky, and R.F. Swarttouw. {\it Hypergeometric orthogonal polynomials and their q-analogues}. Springer, 1-st edition, 2010.

\bibitem{MNR} A.Magnus, F.Ndayiragije and A.Ronveaux. {\it Heun differential equation satisfied by some classical
biorthogonal \\ rational functions} (to be published); https://perso.uclouvain.be/alphonse.magnus/num3/biorthclassCanterb2017.pdf

\bi{NT} K.Nomura and P.Terwilliger, {\it Linear transformations that are tridiagonal with respect to both eigenbases of a Leonard pair}, Lin.Alg.Appl.
{\bf 420} (2007), 198--207.  arXiv:math/0605316.

\bibitem{Ronveaux} A. Ronveaux (Ed.), {\it Heun's Differential Equations}, Oxford University Press, Oxford, 1995.

\bibitem{Ros} H.Rosnegren, {\it An elementary approach to $6j$-symbols (classical, quantum, rational, trigonometric, and elliptic)}, Ramanujan J., {\bf 13} (2007), 131–-166, arXiv:math/0312310. 


\bi{Takemura} K. Takemura {\it On $q$-deformations of Heun equation}, Arxiv:1712.09564.

\bi{Ter} P.Terwilliger, {\it Two linear transformations each tridiagonal with respect to an eigenbasis of the other}, Lin.Alg.Appl. {\bf 330} (2001), 149–-203.

\bi{TVZ} S.Tsujimoto, L.Vinet and A.Zhedanov {\it Dunkl shift operators and Bannai-Ito polynomials}, Advances in Mathematics
{\bf 229}, (2012), 2123-–2158.

\bi{Tur_Heun} A.Turbiner, {\it The Heun operator as a Hamiltonian}, Journal of Physics {\bf A49} (2016), 26LT01. arXiv:1603.02053

\bi{Tur_quasi} A. Turbiner, {\it One-Dimensional Quasi-Exactly Solvable Schr\"odinger Equations}, Physics Reports {\bf 642} (2016) 1--71.  arXiv:1603.02992.

\bi{Zh_AW} A. S. Zhedanov, {\it "Hidden symmetry" of Askey-Wilson polynomials}, Theoret. and Math. Phys. {\bf 89}
(1991), 1146–-1157.

\bibitem{Zhe_Pade} A.Zhedanov, {\it Pad\'e interpolation table and biorthogonal rational functions}, RIMS Proc. 2004.

\end{thebibliography}

\end{document}